\documentclass[11pt,twoside,a4paper]{article}
\usepackage{amsfonts}\usepackage{amsmath}
\usepackage{amsthm}\usepackage{a4}\usepackage{url}
\usepackage{amssymb}\usepackage{graphics}
\numberwithin{equation}{section}
\newtheorem{thm}{Theorem}[section]
\newtheorem{cor}[thm]{Corollary}
\newtheorem{lem}[thm]{Lemma}
\theoremstyle{definition}

\newtheorem{rem}[thm]{Remark}

\begin{document}
\title{A note on ergodic transformations of self-similar Volterra Gaussian processes}\author{C\'eline Jost\footnote{Department of Mathematics and Statistics, P.O.\,Box 68 (Gustaf H\"allstr\"omin katu
2b), 00014 University of Helsinki, Finland. E-mail: celine.jost@iki.fi}}\date{}\maketitle
\begin{abstract}
We derive a class of ergodic transformations of self-similar Gaussian processes that are Volterra, i.e. of type $X_t=\int^t_0 z_X(t,s)dW_s$, $t\in[0,\infty)$, where $z_X$ is a deterministic kernel and $W$ is a standard Brownian motion.
\end{abstract}
\textit{MSC:} 60G15; 60G18; 37A25\\
\\
\textit{Keywords:} Volterra Gaussian process; Self-similar process; Ergodic transformation; Fractional Brownian motion 
\section{Introduction}
Let $(X_t)_{t\in[0,\infty)}$ be a continuous \textit{Volterra} Gaussian process on a complete probability space $(\Omega,\mathcal{F},\mathbb{P})$. This means that \begin{equation}X_t\ =\ \int^{t}_0 z_X(t,s)dW_s,\ a.s.,\ t\in[0,\infty),\label{X=intzdW}\end{equation} where the kernel $z_X\in L^2_{\text{loc}}\left([0,\infty)^2\right)$ is \textit{Volterra}, i.e. $z_X(t,s)=0$, $s\geq t$, and $(W_t)_{t\in[0,\infty)}$ is a standard Brownian motion. Clearly, $X$ is centered and \begin{equation*}R^{X}(s,t)\ :=\ \text{Cov}_{\mathbb{P}}\left(X_s,X_t\right)\ =\ \int^s_0 z_X(s,u)z_X(t,u)du,\ 0\leq s\leq t<\infty.\label{RX=intz_X}\end{equation*}
We assume that $X$ is $\beta$\textit{-self-similar} for some $\beta>0$, i.e. \begin{equation}\left(X_{at}\right)_{t\in[0,\infty)}\ \stackrel{d}{=}\ \left(a^{\beta}X_{t}\right)_{t\in[0,\infty)},\ a>0,\label{selfsimilarfdd}\end{equation}where $\stackrel{d}{=}$ denotes equality of finite-dimensional distributions.
Furthermore, we assume that $z_X$ is \textit{non-degenerate} in the sense that the family $\{z_X(t,\cdot)\,|\,t\in(0,\infty)\}$ is linearly independent and generates a dense subspace of $L^2\left([0,\infty)\right)$. Then \begin{equation}\Gamma_t(X)\ :=\ \overline{\text{span}\{X_s\,|\,s\in[0,t]\}}\ = \ \Gamma_t(W),\  t\in(0,\infty),\label{chsX=chsW}\end{equation} where the closure is in $L^2(\mathbb{P})$, or equivalently, $$\mathbb{F}^X\ =\ \mathbb{F}^W,$$ where $\mathbb{F}^X:=\left(\mathcal{F}^X_{t}\right)_{t\in[0,\infty)}$ denotes the completed natural filtration of $X$.\\
\\
We assume implicitly that $(\Omega,\mathcal{F},\mathbb{P})$ is the \textit{coordinate space of} $X$, which means that $\Omega=
\{\omega:[0,\infty)\to\mathbb{R}\,|\,\omega\text{ is continuous}\}$, $\mathcal{F}=\mathcal{F}^{X}_{\infty}:=\sigma(X_t\,|\,t\in[0,\infty))$ and $\mathbb{P}$ is the probability measure with respect to which the \textit{coordinate process} $X_t(\omega)=\omega(t)$, $\omega\in \Omega$, $t\in[0,\infty)$, is a centered Gaussian process with covariance function $R^X$.  Recall that a measurable map \begin{eqnarray*}\mathcal{Z}\ :\  (\Omega,\mathcal{F},\mathbb{P}) &\to&  (\Omega,\mathcal{F},\mathbb{P})\\
X(\omega)\quad &\mapsto& \mathcal{Z}(X(\omega))
\end{eqnarray*}
 is a \textit{measure-preserving transformation}, or \textit{endomorphism, on} $(\Omega,\mathcal{F},\mathbb{P})$ if $\mathbb{P}^{\mathcal{Z}}=\mathbb{P}$, or equivalently, if $\mathcal{Z}(X)\stackrel{d}{=}X$. If $\mathcal{Z}$ is also bijective and $\mathcal{Z}^{-1}$ is measurable, then it is an \textit{automorphism}. \\
\\
A classical example for a process of the above type is the widely studied \textit{fractional Brownian motion with Hurst index} $H\in(0,1)$, or $H$-\textit{fBm}, denoted by $\bigl(B^H_t\bigr)_{t\in[0,\infty)}$. $H$-fBm is the continuous, centered Gaussian process with covariance function $$R^{B^H}(s,t)\ =\ \frac{1}{2}\left(s^{2H}\ +\ t^{2H}\ -\ |s-t|^{2H}\right),\ s,t\in[0,\infty).$$ For $H=\frac{1}{2}$, fBm is standard Brownian motion. $H$-fBm is $H$-self-similar and has stationary increments. The non-degenerate Volterra kernel is given by
$$z_{B^H}(t,s)\ =\ c(H)(t-s)^{H-\frac{1}{2}}\cdot{}_{2} F_1\left(\frac{1}{2}-H,H-\frac{1}{2},H+\frac{1}{2},1-\frac{t}{s}\right),\ 0<s<t<\infty,$$ where $c(H):=\left(\frac{2H\Gamma\left(\frac{3}{2}-H\right)}{\Gamma\left( H+\frac{1}{2}\right)\Gamma\left(2-2H\right)}\right)^{\frac{1}{2}}$ with $\Gamma$ denoting the gamma function, and ${}_2F_1$ is the Gauss hypergeometric function.
In 2003, Molchan showed that the transformation \begin{equation}\mathcal{Z}_t\left(B^H\right)\ :=\ B^H_t\ -\ 2H\int^t_0 \frac{B^H_s}{s}ds,\ t\in[0,\infty),\label{MolchansfBm}\end{equation} is measure-preserving and satisfies
  \begin{equation}\Gamma_T\left(\mathcal{Z}\bigl(B^H\bigr)\right)\ =\ \Gamma_T\left(Y^H\right),\ T>0,\label{GammaZBH=GammaYH}\end{equation} where \begin{equation}Y^H_t \ :=\ M^H_t\ -\ \frac{t}{T} \xi^H_T,\ t\in[0,T].\label{YH=}\end{equation}
Here, $M^H_t:=\sqrt{2-2H}\int^t_0 s^{\frac{1}{2}-H}dW_s$, $t\in[0,\infty)$, is the \textit{fundamental martingale of $B^H$} and $\xi^{H}_T:= 2H\int^T_0\left(\frac{s}{T}\right)^{2H-1}dM^H_s$. \\
\\
In this work, we present a class of measure-preserving transformations (on the coordinate space) of $X$, which generalizes this result. Moreover, we show that these measure-preserving transformations are ergodic.  
\section{Ergodic transformations}
First, we introduce the class of measure-preserving transformations:
\begin{thm}\label{lemmpt}Let $\alpha>\frac{-1}{2}$. Then the transformation \begin{equation}\mathcal{Z}^{\alpha}_t(X)\ :=\ X_t \ -\ (2\alpha+1)t^{\beta-\alpha-\frac{1}{2}}\int^t_0 s^{\alpha-\beta-\frac{1}{2}}X_s ds,\ t\in[0,\infty),\label{zalpha:=}\end{equation}
is an automorphism on the coordinate space of $X$. The  inverse is given by
\begin{equation}\mathcal{Z}^{\alpha,-1}_t(X)\ =\ X_t\ -\ (2\alpha+1)t^{\alpha+\beta+\frac{1}{2}}\int^{\infty}_t X_s s^{-\beta-\alpha-\frac{3}{2}}ds,\ a.s.,\ t\in[0,\infty).\label{z-alpha:=}\end{equation}
The integrals on the right-hand sides are $L^2(\mathbb{P})$-limits of Riemann sums.
\end{thm}
\begin{proof}First, note that $R^X$ is continuous. Furthermore, by combining H\"older's inequality and (\ref{selfsimilarfdd}), we have that $$\big|R^X(s,t)\big|\ \leq\ \text{E}_{\mathbb{P}}(X_1)^2 s^{\beta}t^{\beta},\ s,t\in(0,\infty).$$
Hence, the double Riemann integrals $$\int^t_0 \int^t_0 (us)^{\alpha-\beta-\frac{1}{2}}R^X(u,s)duds$$ and 
$$\int^{\infty}_t\int^{\infty}_t (us)^{-\beta-\alpha-\frac{3}{2}}R^X(u,s)duds$$ are finite.
Thus, the integrals in (\ref{zalpha:=}) and (\ref{z-alpha:=}) are well-defined (see \cite{Hu}, section 1).\\
Second, we show that $\mathcal{Z}^{\alpha}$ is a measure-preserving transformation. Let $Y_t:=\exp(-\beta t)X_{\exp(t)}$ and $Y^{\alpha}_t:=\exp(-\beta t)\mathcal{Z}^{\alpha}_{\exp(t)}(X)$, $t\in\mathbb{R}$, denote the Lamperti transforms of $X$ and $\mathcal{Z}^{\alpha}(X)$, respectively. Hence, the process $(Y_t)_{t\in\mathbb{R}}$ is stationary. By substituting $v:=\ln(s)$, we obtain that
\begin{eqnarray*}
Y^{\alpha}_t 
&=& \exp(-\beta t)\left(X_{\exp(t)}-(2\alpha+1)\exp\left(\left(\beta-\alpha-\frac{1}{2}\right)t\right)\int^{\exp (t)}_0 s^{\alpha-\beta-\frac{1}{2}}X_s ds\right)\\
&=&Y_t- (2\alpha +1)\exp\left(\left(-\alpha-\frac{1}{2}\right) t\right)\int^{t}_{-\infty} \exp\left(v\left(\alpha-\beta+\frac{1}{2}\right)\right)X_{\exp (v)} dv\\
&=&Y_t-(2\alpha+1)\exp\left(\left(-\alpha-\frac{1}{2}\right)t\right)\int^t_{-\infty}\exp\left(v\left(\alpha+\frac{1}{2}\right)\right)Y_v dv\\
&=&\int^{\infty}_{-\infty}h^{\alpha}(t-v) Y_v dv,\ a.s.,\ t\in\mathbb{R},\end{eqnarray*}
where $$h^{\alpha}(x)\ :=\ \delta_{0}(x) \ -\  (2\alpha+1)1_{(0,\infty)}(x) \exp\left(-\left(\alpha+\frac{1}{2}\right)x\right),\ x\in\mathbb{R}.$$ Thus, $Y^{\alpha}$ is a linear, non-anticipative, time-invariant transformation of $Y$.
The spectral distribution function of $Y^{\alpha}$ is given by (see \cite{Ya}, p. 151) $$dF^{\alpha}(\lambda)\ =\ \big|H^{\alpha}(\lambda)\big|^2 dF(\lambda),\ \lambda\in\mathbb{R},$$ where $H^{\alpha}(\lambda):=\int_{\mathbb{R}}\exp(-i\lambda x)h^{\alpha}(x)dx$ denotes the Fourier transform of $h^{\alpha}$ and $F$ is the spectral distribution function of $Y$.  
We have that (see \cite{Er1}, p. 14 and p. 72) $$\big|H^{\alpha}(\lambda)\big|\ =\ \Bigg|1 - (2\alpha+1) \left(\frac{\alpha+\frac{1}{2}-i\lambda}{\left(\frac{1}{2}+\alpha\right)^2+\lambda^2}\right)\Bigg|\ =\
\bigg|\frac{-\frac{1}{2}-\alpha+i\lambda}{\frac{1}{2}+\alpha+i\lambda}\bigg|\ =\ 1,\ \lambda\in\mathbb{R},$$
i.e. $F^{\alpha}\equiv F$. It follows from this that $\left(Y^{\alpha}_t\right)_{t\in\mathbb{R}}\stackrel{d}{=}(Y_t)_{t\in\mathbb{R}}$, or equivalently, $\left(\mathcal{Z}^{\alpha}_t(X)\right)_{t\in[0,\infty)}\stackrel{d}{=}(X_t)_{t\in[0,\infty)}$.\\
Third, by splitting integrals and using Fubini's theorem, we obtain that $$\mathcal{Z}^{\alpha,-1}_t\left(\mathcal{Z}^{\alpha}( X)\right)\ =\ X_t\ =\ \mathcal{Z}^{\alpha}_t\left(\mathcal{Z}^{\alpha,-1} (X)\right),\ a.s.,\ t\in[0,\infty).\qedhere$$\end{proof}
\begin{rem}Theorem \ref{lemmpt} generalizes (\ref{MolchansfBm}). In fact, $\mathcal{Z}^{H-\frac{1}{2}}\left(B^H\right)=\mathcal{Z}\left(B^H\right)$, $H\in(0,1)$.\end{rem} 
\begin{rem}Theorem \ref{lemmpt} holds true for general continuous centered $\beta$-self-similar Gaussian processes.\end{rem}
Next, we present two auxiliary lemmas concerning the structure of $z_X$:
\begin{lem}\label{sskernel3points}Let $(X_t)_{t\in[0,\infty)}$ be a Volterra Gaussian process with a non-degenerate Volterra kernel $z_X$. Then the following are equivalent:\\
1. $X$ is $\beta$-self-similar, i.e. 
   \begin{equation*}\int^s_0 z_X(at,au)z_X(as,au)du\  =\ a^{2\beta-1}\int^s_0 z_X(t,u)z_X(s,u)du,\ 0<s\leq t<\infty,\ a>0.\label{equivcond0}\end{equation*}
2. It holds that \begin{equation*}z_X(at,as)\ =\ a^{\beta-\frac{1}{2}}z_X(t,s),\ 0<s<t<\infty,\ a>0.\label{equivcond1}\end{equation*}
3. There exists $F_X\in L^2\left((0,1),(1-x)^{2\beta-1}dx\right)$ such that \begin{equation*}z_X(t,s)\ =\ (t-s)^{\beta-\frac{1}{2}}F_X\left(\frac{s}{t}\right),\ 0<s<t<\infty.\label{equivcond2}\end{equation*}
\end{lem}
\begin{proof}$1\Rightarrow 2$: For $a>0$, let $z_{Y(a)}(t,s):=a^{\frac{1}{2}-\beta}z_X(at,as)$, $0<s<t<\infty$, and let $Y_t(a):=\int^t_0 z_{Y(a)}(t,s)dW_s$, $t\in[0,\infty)$. Clearly, $z_{Y(a)}$ is non-degenerate. From (\ref{chsX=chsW}), we obtain that $\Gamma_t(Y(a))=\Gamma_t(W),$ $t\in(0,\infty)$. From part 1, it follows that $X\stackrel{d}{=}Y(a)$. Thus, the process $W'_t:=\int^t_0 z^{\ast}_X(t,s)dY_s(a)$, $t\in[0,\infty)$, where $z_X^{\ast}$ is the reciprocal of $z_X$ and the integral is an abstract Wiener integral, is a standard Brownian motion with $\Gamma_t(W')=\Gamma_t(Y(a))$, $t\in(0,\infty)$. Hence, $\Gamma_t(W)=\Gamma_t(W')$, i.e. $W$ and $W'$ are indistinguishable. Therefore, $Y_t(a)=\int^t_0 z_X(t,s)dW_s$, a.s., $t\in(0,\infty)$, i.e. $Y_t(a)=X_t$, a.s., $t\in[0,\infty)$. In particular, $0=\text{E}_{\mathbb{P}}(Y_t(a)-X_t)^2 =\int^{t}_0 \bigl(z_{Y(a)}(t,s)-z_X(t,s)\bigr)^2 ds$, $t\in(0,\infty)$. Thus, $z_X(t,\cdot)\equiv z_{Y(a)}(t,\cdot)$, $t\in(0,\infty)$.\\
$2\Rightarrow 3$: Let $G_X(t,s):=(t-s)^{\frac{1}{2}-\beta}z_X(t,s)$, $0<s<t<\infty$. From part 2, it follows that $G_X(at,as)=G_X(t,s),$ $0<s<t<\infty$, $a>0$. Hence, for every $(t,s)$, $s<t$, the function $G_X$ is constant on the line $\{(at,as)\,|\,a\in(0,\infty)\}$, which depends only on the slope $\frac{s}{t}$. Thus, $G_X(t,s)=F_X\bigl(\frac{s}{t}\bigr)$, $0<s<t<\infty$, for some $F_{X}\in L^2\left((0,1),(1-x)^{2\beta-1}dx\right)$.\\
$3\Rightarrow 1$: This is trivial.  \end{proof}
\begin{lem}\label{kernelstructure}Let $\alpha>\frac{-1}{2}$. Then we have that \begin{equation*} t^{\beta-\alpha-\frac{1}{2}}\int^t_s  u^{\alpha-\beta-\frac{1}{2}}z_X(u,s)du\ =\ s^{\alpha}\int^t_s z_X(t,u)u^{-\alpha-1}du,\ 0<s<t<\infty.\label{kernelidentity}\end{equation*}\end{lem}
\begin{proof}From Lemma \ref{sskernel3points}, it follows that the Volterra kernel of $X$ can be written as \begin{equation*}z_X(t,s)\ =\ (t-s)^{\beta-\frac{1}{2}}F_X\left(\frac{s}{t}\right),\ 0<s<t<\infty,\label{kernelfactorization}\end{equation*}for some function $F_X$.
By substituting first $x:=\frac{s}{u}$ and then $v:=tx$, we obtain that \begin{eqnarray*}t^{\beta-\alpha-\frac{1}{2}}\int^t_s u^{\alpha-\beta-\frac{1}{2}}z_X(u,s)du&=&t^{\beta-\alpha-\frac{1}{2}}\int^t_s u^{\alpha-\beta-\frac{1}{2}}(u-s)^{\beta-\frac{1}{2}}F_X\left(\frac{s}{u}\right)du\\
&=&t^{\beta-\alpha-\frac{1}{2}}s^{\alpha}\int^1_{\frac{s}{t}}x^{-\alpha-1}(1-x)^{\beta-\frac{1}{2}}F_X(x)dx\\
&=&s^{\alpha}\int^t_s (t-v)^{\beta-\frac{1}{2}}F_X\left(\frac{v}{t}\right)v^{-\alpha-1}dv\\
&=&s^{\alpha}\int^t_s z_X(t,v)v^{-\alpha-1}dv.\qedhere\end{eqnarray*} 
\end{proof}
The next lemma is the key result for deriving the ergodicity of the measure-preserving transformations:
\begin{lem}\label{ZalphaX=intZalphaW}Let $\alpha>\frac{-1}{2}$. Then $$\mathcal{Z}^{\alpha}_t(X)\ =\ \int^t_0 z_X(t,s)d\mathcal{Z}^{\alpha}_s(W),\ a.s., \ t\in[0,\infty).$$
\end{lem}
\begin{proof}By combining (\ref{X=intzdW}) and the stochastic Fubini theorem, using Lemma \ref{kernelstructure}, again the stochastic Fubini theorem, and finally using partial integration, we obtain that 
 \begin{eqnarray*}\mathcal{Z}^{\alpha}_t(X)&=& X_t\ -\ (2\alpha+1)\int^t_0\left(t^{\beta-\alpha-\frac{1}{2}}\int^t_u s^{\alpha-\beta-\frac{1}{2}}z_X(s,u)ds\right)dW_u \\
&=&   X_t\ -\ (2\alpha+1)\int^t_0\left(u^{\alpha}\int^t_u z_X(t,s)s^{-\alpha-1}ds\right)dW_u \\
&=&X_t\ -\ (2\alpha+1)\int^t_0 z_X(t,s)s^{-\alpha-1}\int^s_0 u^{\alpha}dW_u ds\\
&=&X_t\ -\ (2\alpha+1)\int^t_0 z_X(t,s)\left((-\alpha)s^{-\alpha-1}\int^s_0 u^{\alpha-1}W_ududs+s^{-1}W_sds\right)\\
&=&\int^t_0 z_X(t,s)d\mathcal{Z}^{\alpha}_s(W),\ a.s.,\ t\in(0,\infty).\qedhere\end{eqnarray*}
\end{proof}
In the following, let $\mathcal{Z}^{\alpha,n}:=\left(\mathcal{Z}^{\alpha}\right)^n$ denote the $n$-th iterate of $\mathcal{Z}^{\alpha}$, $n\in\mathbb{Z}$. Also, let  $$\Gamma_{\infty}(X)\ :=\ \overline{\mathrm{span}\{X_t\ |\ t\in[0,\infty)\}}.$$
For $\alpha>-\frac{1}{2}$, let \begin{equation*}N_t^{\alpha}\ :=\ \int^t_0 s^{\alpha}dW_s,\ t\in[0,\infty).\label{defNalpha}\end{equation*}  
Clearly, $N^{\alpha}$ is an $\left(\alpha+\frac{1}{2}\right)$-self-similar $\mathbb{F}^X$-martingale. From Lemma \ref{ZalphaX=intZalphaW}, it follows that
\begin{eqnarray}\mathcal{Z}^{\alpha}_t(X)&=&\int^t_0 z_X(t,s)s^{-\alpha}d\mathcal{Z}^{\alpha}_s\left(N^{\alpha}\right),\ a.s.,\ t\in(0,\infty).\label{ZTX=intZTN}
\end{eqnarray} 
The next lemma is an auxiliary result, which was obtained in \cite{Jo}, section 3.2 and Theorem 5.2. (The automorphism $\mathcal{Z}^{\alpha}$ on the coordinate space of $N^{\alpha}$ here corresponds to the (ergodic) automorphism $\mathcal{T}^{(1)}$ on the coordinate space of the martingale $M$ with $M:=N^{\alpha}$ in \cite{Jo}.)
\begin{lem}\label{lemma27Na}Let $\alpha>\frac{-1}{2}$ and $T>0$. \\
1. It holds that \begin{equation*}\Gamma_T\left(\mathcal{Z}^{\alpha}\left(N^{\alpha}\right)\right)\ =\ \Gamma_T\left(N^{\alpha,T}\right),\label{partfrom3}\end{equation*} 
where $N^{\alpha,T}_t:=N_t^{\alpha}-\left(\frac{t}{T}\right)^{2\alpha+1}N^{\alpha}_T$, $t\in[0,T]$, is a \textit{bridge of} $N^{\alpha}$, i.e. a process satisfying $\mathrm{Law}_{\mathbb{P}}\left(N^{\alpha,T}\right)=\mathrm{Law}_{\mathbb{P}}\left(N^{\alpha}\, |\, N^{\alpha}_T=0\right)$, and $\Gamma_T\left(N^{\alpha,T}\right):=\overline{\mathrm{span}\big\{N^{\alpha,T}_t\,|\,t\in[0,T]\big\}}$.\\
2. We have that
 \begin{equation}\Gamma_{T}\left(N^{\alpha}\right)\
 =\ \overline{\perp_{n\in\mathbb{N}_0}\mathrm{span}\!\left\{\mathcal{Z}^{\alpha,n}_T\bigl(N^{\alpha}\bigr)\right\}}\label{GammaT(N)}\end{equation}
 and \begin{equation}\Gamma_{\infty}\left(N^{\alpha}\right)\
 =\ \overline{\perp_{n\in\mathbb{Z}}\mathrm{span}\!\left\{\mathcal{Z}^{\alpha,n}_T\bigl(N^{\alpha}\bigr)\right\}}.\label{Gammainfty(N)} \end{equation}
Here, $\perp$ denotes the orthogonal direct sum.
\end{lem}
By combining (\ref{ZTX=intZTN}) and part 1 of Lemma \ref{lemma27Na}, we obtain the following: 
\begin{lem}\label{lemma27X}Let $\alpha>\frac{-1}{2}$ and $T>0$. Then $$\Gamma_T\left(\mathcal{Z}^{\alpha}\left(X\right)\right)\ =\ \Gamma_T\left(N^{\alpha,T}\right).$$
\end{lem}
\begin{rem}Lemma \ref{lemma27X} is a generalization of identity (\ref{GammaZBH=GammaYH}). Indeed, we have that $Y^H_t =\sqrt{2-2H}\int^t_0 s^{1-2H}dN^{H-\frac{1}{2},T}_s$, a.s., $t\in[0,T]$, where $Y^H$ is the process defined in (\ref{YH=}). $Y^H$ is a bridge (of some process) if and only if $H=\frac{1}{2}$.\end{rem}
The following generalizes part 2 of Lemma \ref{lemma27Na}:
\begin{lem}\label{chaos}Let $\alpha>\frac{-1}{2}$ and $T>0$. Then we have that
$$\Gamma_T\left(X\right)\ =\ \overline{\oplus_{n\in\mathbb{N}_0}\,\mathrm{span}\left\{\mathcal{Z}^{\alpha,n}_T(X)\right\}}$$ 
and  $$\Gamma_{\infty}\left(X\right)\ =\ \overline{\oplus_{n\in\mathbb{Z}}\,\mathrm{span}\left\{\mathcal{Z}^{\alpha,n}_T(X)\right\}}.$$\end{lem}
\begin{proof}We assume that $X\neq N^{\alpha}$.
By iterating (\ref{ZTX=intZTN}) and (\ref{GammaT(N)}), we obtain that \begin{equation*}\mathcal{Z}_T^{\alpha,n}(X)\ \in\ \Gamma_{T}\left(\mathcal{Z}^{\alpha,n}\left(X\right)\right)
 =\ \Gamma_{T}\left(\mathcal{Z}^{\alpha,n}\bigl(N^{\alpha}\bigr)\right)\ =\ \overline{\perp_{i\geq n}\text{span}\left\{\mathcal{Z}^{\alpha,i}_T\bigl(N^{\alpha}\bigr)\right\}},\ n\in\mathbb{Z}\label{GammaTN}.\end{equation*}
 Moreover, $X_T\not\perp N^{\alpha}_T$, hence
$\mathcal{Z}^{\alpha,n}_T(X)\not\perp\mathcal{Z}^{\alpha,n}_T\left( N^{\alpha}\right)$, $n\in\mathbb{Z}$, and therefore, 
$$\mathcal{Z}^{\alpha,n}_T(X)\ \not\in\ \overline{\perp_{i\geq n+1}\text{span}\left\{\mathcal{Z}^{\alpha,i}_T\left(N^{\alpha}\right)\right\}},\ n\in\mathbb{Z}.$$  
From (\ref{GammaT(N)}) and (\ref{Gammainfty(N)}), it follows that the systems $\left\{\mathcal{Z}^{\alpha,n}_T(X)\right\}_{n\in\mathbb{N}_0}$ and $\left\{\mathcal{Z}^{\alpha,n}_T(X)\right\}_{n\in\mathbb{Z}}$ are free and complete in $\Gamma_T(X)$ and $\Gamma_{\infty}(X)$, respectively.
\end{proof}
\begin{rem}
The process $X$ is an $\mathbb{F}^X$-Markov process if and only if there exists $\alpha>\frac{-1}{2}$ and a constant $c(X)$, such that \begin{equation}X_t\ =\ c(X)\cdot t^{\beta-\frac{1}{2}-\alpha}\int^t_0 s^{\alpha}dW_s,\ a.s.,\ t\in(0,\infty).\label{markovrepresentation}\end{equation}The free complete system $\left\{\mathcal{Z}^{\alpha,n}_T(X)\right\}_{n\in\mathbb{Z}}$ is orthogonal if and only if (\ref{markovrepresentation}) is satisfied.
\end{rem}
From Lemma \ref{chaos}, we obtain the following:
\begin{cor}\label{sigmaalgX}Let $\alpha>\frac{-1}{2}$ and $T>0$. Then
$$\mathcal{F}^X_T\ =\ \bigvee_{n\in\mathbb{N}_0}\sigma\bigl(\mathcal{Z}^{\alpha,n}_T(X)\bigr).$$
Furthermore, \begin{equation*}\mathcal{F}\ =\ \mathcal{F}^X_{\infty}\ =\ \bigvee_{n\in\mathbb{Z}}\sigma\bigl(\mathcal{Z}^{\alpha,n}_T(X)\bigr).\label{F=VninZ}\end{equation*}
\end{cor}
Recall that an automorphism $\mathcal{Z}$ is a \textit{Kolmogorov automorphism}, if there exists a $\sigma$-algebra $\mathcal{A}\subseteq\mathcal{F}$, such that $\mathcal{Z}^{-1}\mathcal{A}\subseteq\mathcal{A}$, $\vee_{m\in\mathbb{Z}}\mathcal{Z}^m \mathcal{A}=\mathcal{F}$ and $\cap_{m\in\mathbb{N}_0}\mathcal{Z}^{-m}\mathcal{A}=\{\Omega,\emptyset\}$.
A Kolmogorov automorphism is strongly mixing and hence ergodic (see \cite{Pe}, Propositions 5.11 and 5.9 on p. 63 and p. 62). The ergodicity of $\mathcal{Z}^{\alpha}$ is hence a consequence of the following: 
\begin{thm}Let $\alpha>\frac{-1}{2}$ and $T>0$. The automorphisms $\mathcal{Z}^{\alpha}$ and $\mathcal{Z}^{\alpha,-1}$ are Kolmogorov automorphisms with $\mathcal{A}=\vee_{n\in\mathbb{-N}}\sigma\bigl(\mathcal{Z}^{\alpha,n}_T(X)\bigr)$ and $\mathcal{A}=\mathcal{F}^X_T$, respectively.\end{thm}
\begin{proof}
$\mathcal{Z}^{\alpha}$ is a Kolmogorov automorphism with  $\mathcal{A}=\vee_{n\in -\mathbb{N}}\sigma\bigl(\mathcal{Z}^{\alpha,n}_T(X)\bigr)$:\\ 
First, $\mathcal{Z}^{\alpha,-1}\mathcal{A}=\vee_{n\in-\mathbb{N}}\sigma\bigl(\mathcal{Z}^{\alpha,n-1}_T(X)\bigr)\subseteq\mathcal{A}$.\\
Second, $\vee_{m\in\mathbb{Z}}\mathcal{Z}^{\alpha,m}\mathcal{A} = \vee_{m\in\mathbb{Z}}\vee_{n\in-\mathbb{N}}\sigma\bigl(\mathcal{Z}^{\alpha,m+n}_T(X)\bigr)=\mathcal{F}$.\\
Third, let $\{Y_n\}_{n\in-\mathbb{N}}$ denote the Hilbert basis of $\overline{\oplus_{n\in-\mathbb{N}}\,\text{span}\left\{\mathcal{Z}^{\alpha,n}_T(X)\right\}}$ which is obtained from $\left\{\mathcal{Z}^{\alpha,n}_T(X)\right\}_{n\in-\mathbb{N}}$ via Gram-Schmidt orthonormalization. By using Kolmogorov's zero-one law (see \cite{Sh}, p. 381), we obtain that
 $\cap_{m\in \mathbb{N}_0}\mathcal{Z}^{\alpha,-m}\mathcal{A}=\cap_{m\in\mathbb{N}_0}\left(\vee_{n\leq -m-1}\sigma\bigl(\mathcal{Z}^{\alpha,n}_T(X)\bigr)\right)=\cap_{m\in\mathbb{N}_0}\left(\vee_{n\leq -m-1}\sigma\bigl(Y_n\bigr)\right)
=\{\Omega,\emptyset\}$.\\
Similarly, one shows that $\mathcal{Z}^{\alpha,-1}$ is a Kolmogorov automorphism with $\mathcal{A}=\mathcal{F}^X_T$.
\end{proof}
\textbf{Acknowledgements.} Thanks are due to my supervisor Esko Valkeila and to Ilkka Norros for helpful comments. I am indebted to the Finnish Graduate School in Stochastics (FGSS) and the Finnish Academy of Science and Letters, Vilho, Yrj\"o and Kalle V\"ais\"ala Foundation for financial support.

\end{document}